# Twofold fast summation

Evgeny Latkin
2013 Dec 31

Debugging accumulation of floating-point errors is hard; ideally, computer should track it automatically. Here we consider twofold approximation of exact real with *value + error* pair of floating-point numbers. Normally, *value + error* sum is more accurate than *value* alone, so *error* can estimate deviation between *value* and its exact target. Fast summation algorithm, that provides twofold sum of $\sum x_n$ or dot product $\sum x_n y_n$, can be same fast as direct summation if leveraging processor underused potential. This way, we hit three goals: improve precision, track inaccuracy, and do this with little if any loss in performance.

Contents

Document history:
30-Nov-2013: first version
31-Dec-2013: bug fix: misprint in Knuth theorem thus wrong code for rigorous twofold summation

## Overview

Here I would seek for a reasonably good compromise for three goals:

- Improve precision of floating-point calculations
- Automate tracking accumulation of rounding errors
- Do this extra job with little if any loss in performance

Accumulation of floating-point inaccuracy is difficult to debug; so ideally, computer could automatically track rounding errors and adapt precision if necessary or at least signal if precision appears not enough. People tried many methods here [1], of which I would mention intervals [2] and just increasing floating-point precision, maybe with double-double or similar techniques [3], [4].

Here I consider simple "twofold" approach representing an exact target value with *value* + *error* sum of floating-point numbers, presumably of IEEE-754 single or double precision. Such approximation mixes interval and double-double techniques. Normally, *value* + *error* sum is more accurate than *value* alone, so *error* can reasonably estimate the interval between *value* and its exact counterpart.

Unlike intervals, twofold approximation cannot guarantee if exact solution lays within the *error* interval, and can lead to completely wrong perception of what the error is. But I think more probably, computing with twofolds would signal if accumulated error gets too large, so you can interrupt and consider using higher precision. The benefit of such less rigorous approach is that twofolds must not suffer redundant widening like intervals striving to cover all possible errors including the worst case. Twofolds rather address average-case like ordinary floating-point numbers, and ultimately can be exact.

Unlike double-double approach, with twofold we accent on tracking inaccuracy rather than reducing it. This way, twofold is normally less accurate than could be with double-double renormalization. In turn, computing with twofolds is simpler so potentially can demonstrate higher performance.

According to performance testing with my commodity laptop (Ivy Bridge processor), twofold summation for $\sum x_n$ or $\sum x_n y_n$ can be very fast. With best of my hand-made AVX vectorization, twofold summation with double precision yields 3 gigaflop on one CPU core at 3 GHz, which is only 4 times slower than this processor theoretical peak 12 gigaflop. This is 100 times better than 30 megaflop with g++ __float128.

If we take into account the cost of memory reading, twofold loses only 2 times if summing a larger array that cannot fit into fastest L1 cache, and twofold does not lose at all if array can't fit into last-level (L3) cache. If at least part of data is not pre-loaded into cache, CPU capacity is underused while fetching data from memory, and twofold algorithm can leverage this resource.

Average-case inaccuracy estimate for sum of $N$ numbers is $O(N\varepsilon^2)$ for twofold *value* + *error* result, which must be better than $O(\sqrt{N}\varepsilon)$ for direct sum if $N$ is not too large. As we want *error* to estimate deviation from exact result, we need *value* + *error* be much more accurate than *value* alone, that is $O(N\varepsilon^2)$ be much less than $O(\sqrt{N}\varepsilon)$, so $N$ be much less than $1/\varepsilon^2$. I guess, twofold must work fine if $N < 1/\varepsilon$, that is $N$ must not exceed $2^{24} \approx 10^7$ for single precision and $2^{53} \approx 10^{16}$ for double.

## Performance

Key advantage of twofold fast summation is high performance, so let me start with performance results. You can find the test sources and the testing logs at my Web page [7]:

Table 1. Megaflops in various methods of summation and dot product

|  | float | | | double | | |
|---|---|---|---|---|---|---|
|  | small | medium | large | small | medium | large |
| sum | 1141.9 | 1008.73 | 996.973 | 1041.46 | 1012.92 | 982.468 |
| dot | 1090.09 | 1008.73 | 931.707 | 1029.39 | 1012.92 | 870.616 |
| sumtf | 622.326 | 603.757 | 603.319 | 632.612 | 612.368 | 603.757 |
| dottf | 635.437 | 603.319 | 592.445 | 645.769 | 611.399 | 545.601 |
| sumt | 381.682 | 378.978 | 375.008 | 381.682 | 385.235 | 373.293 |
| dott | 381.682 | 379.363 | 367.414 | 381.682 | 381.682 | 355.919 |
| sumk | 272.63 | 252.071 | 252.707 | 272.8 | 255.853 | 251.658 |
| dotk | 272.63 | 252.707 | 248.513 | 272.63 | 255.853 | 244.599 |
| sumd | 1073.74 | 1012.92 | 1002.36 | 29.3601 | - | - |
| dotd | 925.632 | 955.253 | 875.678 | 29.2173 | - | - |
| usum | 2862.61 | 2986.34 | 2763.25 | 2921.33 | 2943.44 | 1886.61 |
| udot | 2477.19 | 2462.06 | 1817.73 | 2440.04 | 1596.21 | 1008.73 |
| usumtf | 735.861 | 759.02 | 745.538 | 738.769 | 745.538 | 738.007 |
| udottf | 734.52 | 747.635 | 693.626 | 745.538 | 756.891 | 670.663 |
| usumt | 431.767 | 428.868 | 428.868 | 425.185 | 428.868 | 425.185 |
| udott | 426.25 | 428.868 | 411.793 | 425.185 | 428.868 | 393.892 |
| usumk | 747.635 | 747.635 | 745.538 | 747.635 | 750.398 | 735.861 |
| udotk | 745.538 | 750.398 | 703.867 | 747.635 | 747.635 | 681.737 |
| usumd | 2031.13 | 2232.42 | 2135.99 | - | - | - |
| udotd | 1275.38 | 1282.5 | 1160.47 | - | - | - |
| vsum | 4658.45 | 7945.85 | 4913.63 | 3120.56 | 3979.63 | 2218.47 |
| vdot | 4169.14 | 4506.79 | 2367.68 | 2942.21 | 2147.48 | 1141.9 |
| vsumd | 4366.66 | 5008.6 | 3939.74 | - | - | - |
| vdotd | 3545.03 | 3870.24 | 2182.4 | - | - | - |
| asum | 8103.4 | 8010.51 | 4958.2 | 4043.31 | 4043.31 | 2200.21 |
| adot | 8072.99 | 4364.8 | 2330.5 | 3435.13 | 2147.48 | 1126.7 |
| asumtf | 4902.09 | 4816.63 | 4228.4 | 2423.28 | 2412.9 | 2028.96 |
| adottf | 4836.25 | 4224.11 | 2253.39 | 2423.28 | 2080.9 | 1112.65 |
| asumt | 3054.5 | 3034.9 | 2862.61 | 1513.37 | 1513.37 | 1448.84 |
| adott | 3034.9 | 2986.34 | 2115.75 | 1521.44 | 1505.2 | 1029.39 |
| asumk | 2022.11 | 2011.33 | 2020.21 | 1008.73 | 1008.73 | 997.963 |
| adotk | 2028.96 | 2031.13 | 1804.89 | 1008.73 | 1011.06 | 902.824 |
| asumd | 4043.31 | 4035.97 | 3638.09 | - | - | - |
| adotd | 3384.37 | 3200.25 | 2115.75 | - | - | - |
| hsum | 23893.9 | 8862.55 | 5019.95 | 11923.9 | 4436.95 | 2138.39 |
| hdot | 12242 | 4568.72 | 2367.68 | 6072.16 | 2216.19 | 1129.08 |
| hsumtf | 5964.3 | 5981.08 | 4436.53 | 2943.44 | 2986.34 | 2080.9 |
| hdottf | 5981.08 | 4043.31 | 2253.39 | 3006.4 | 1997.03 | 1112.65 |
| hsumt | 3319.5 | 3435.13 | 3125.81 | 1720.71 | 1720.71 | 1584.04 |
| hdott | 3435.13 | 3322.44 | 2113.67 | 1720.71 | 1637.23 | 1041 |
| hsumk | 5981.08 | 5981.08 | 4510.43 | 3034.9 | 3014.79 | 2080.9 |
| hdotk | 5981.08 | 4432.37 | 2291.87 | 2986.34 | 2200.21 | 1126.7 |
| hsumd | 5726.27 | 6078.33 | 4432.37 | - | - | - |
| hdotd | 3469.8 | 3410.92 | 2115.75 | - | - | - |
| psum | 23926.4 | - | - | 12049.3 | - | - |
| pdot | 24052.2 | - | - | 12049.3 | - | - |
| psumtf | 5981.08 | - | - | 3031.82 | - | - |
| pdottf | 6072.16 | - | - | 3036.08 | - | - |
| psumt | 3435.13 | - | - | 1720.71 | - | - |
| pdott | 3435.13 | - | - | 1709.18 | - | - |
| psumk | 5981.08 | - | - | 2986.34 | - | - |
| pdotk | 5981.08 | - | - | 3006.4 | - | - |
| psumd | 11923.9 | - | - | - | - | - |
| pdotd | 12049.3 | - | - | - | - | - |
| read1 | 24052.2 | 9013.57 | 5070.82 | 12109 | 4427.09 | 2147.48 |
| read2 | 12109 | 4836.25 | 2316.14 | 6078.33 | 2351.53 | 1112.65 |

For testing, I used my inexpensive Lenovo V580c laptop built with Intel Core i5-3210M processor with nominal frequency at 2.5 GHz and maximal at 3.1 GHz. Memory was two 4GB banks of DDR3-PC12800, so maximal bandwidth was 25.6 gigabytes per second, same as for the processor. I tested single thread to measure performance per CPU core; so processor actually worked at nearly maximal frequency.

I tested with GNU C++ 4.8.1 (Cygwin) and Microsoft C++ 18.0 (Visual Studio Express 2013), and Intel C++ 14.0.1 (Composer XE 2013 SP1). The GNU and Intel compilers support quad precision types __float128 and _Quad, which I used for checking accuracy. This table shows results with GNU compiler; results with Microsoft compiler are essentially same, and results with Intel compiler are higher in some cases.

This table summarizes megaflops I observed in summation and dot product, counting only summation operations. Modern processor can multiply $x_n y_n$ in parallel with summation, so here we count both $\sum x_n$ and $\sum x_n y_n$ like same $N$ operations. We also ignore any additional operations twofold and Kahan summation methods do for assessing or compensating inaccuracy. This way we compare performance of various methods against simple direct summation.

To track dependency on memory bandwidth, I tested with small, medium-size, and large data arrays, which fit to processor fastest L1 cache, fit to last-level (L3) cache, and do not fit to processor cache.

Every function tested with single and double precisions, float and double types of C/C++. For testing with higher-precision accumulator, I used double accumulator for summation of float data, and used quad accumulator for double-precision data (__float128 with GNU C++, and _Quad for Intel C++).

Same functions compared with and without vectoring for AVX. For no-vectoring variant, I used the best of compiler optimization with strict floating-point math. Note, that we cannot allow fast-math compiler optimizations for twofold and Kahan summation, as these methods essentially base on the tricks with correctly rounding the results of the specific math expressions. Fast-math optimization can damage it.

For vectored variant, I have made more-or-less tricky hand-made vectoring with AVX intrinsic functions supplied with the compilers. The reason for optimizing manually is that compiler's optimized code gives only around 1/3 of processor's peak performance. My manual optimization yields higher performance close to processor theoretical peak, up to 24 gigaflops for float and 12 gigaflops for double.

The functions tested here are the following:

- sum, dot: simple direct summation and dot product, compiled for strict math, no AVX
- sumd, dotd: same but with higher-precision accumulator, strict math, no AVX vectoring

- sumtf, dottf: twofold fast summation, strict floating-point math, no AVX vectorization
- sumt, dott: twofold rigorous summation, strict floating-point math, no AVX vectorization
- sumk, dotk: compensatory summation (aka: Kahan summation), strict math, no AVX
- usum, udot, usumd, udotd, usumtf, udottf, usumt, udott, usumk, udotk: strict math, no AVX, but trickier C/C++ coding for peak performance for scalar computing (without SIMD vectorization)
- vsum, vdot, vsumd, vdotd: simple functions vectorized for AVX with best compiler optimization
- asum, adot, asumd, adotd, asumtf, adottf, asumt, adott, asumk, adotk: naïve AVX vectoring
- hsum, hdot, hsumd, hdotd, hsumtf, hdottf, hsumt, hdott, hsumk, hdotk: better AVX vectoring, for performance close to theoretical peak, appears much faster than naïve AVX vectoring
- psum, pdot, psumd, pdotd, psumtk, pdottk, psumt, pdott, psumk, pdotk: same as hsum, … but without actually reading from memory, I use it for checking if can get close to processor peak
- read1, read2: reading memory without actually doing any computations, test memory impact

Let us interpret these performance results. Please note that my test is not designed for exact measuring the execution time, it rather assumes to demonstrate how the concept works in overall:

## Memory reading

Functions read1 and read2 test AVX reading from memory with one and two arrays of data, simulating reading for summation and for dot product correspondingly. For small array, my test shows nearly one reading per CPU tick, which at 3 GHz would give reading 24 (=3×8) billion float numbers or 12 (=3×4) billions of double-precision numbers per second. For summation, this can feed up to 24 gigaflops for float and 12 gigaflop for double, which is the processor peak. For dot product, such reading bandwidth would feed twice fewer gigaflops, because dot product needs twice more data per summation.

Such 2x loss in read2 test looks unexpected for me; maybe coding in assembler could give better results. My reading test and manually optimized functions are written with AVX intrinsic functions supplied with C/C++ compiler, which might imply some limitations.

Reading float numbers through one channel (for summation) from large array that cannot fit into cache shows 5 billion per second, which would be 20 gigabytes per second, fairly close to maximal 25.6 GB per second for this processor and memory. Reading from two large arrays gives twice less floats per second. For double precision, reading shows similar pattern, though at 20% fewer performance for large arrays (such 20% difference looks unexpected for me).

## Purely computing (w/o memory reading)

Consider manually vectorized function names starting with "p" like psum, pdot, etc. These p-functions purely compute but do not actually read from memory (they operate with the minimal piece of data in processor registers). This way, we can estimate quality of optimization: if encoded optimally, p-functions must demonstrate highest performance close to processor theoretical peak. And they actually do:

Theoretical peak for psum, pdot is doing one AVX summation per processor tick, which at 3 GHz would do 24 (=3×8) gigaflops for float and 12 (=3×4) gigaflop for double type. Note, we count only summations, and do not count multiplications $x_n y_n$ for dot product. Look at the table for psum, pdot results; they look pretty close to the peak, modulo some inaccurate measuring in this test.

Twofold fast summation spends 3 extra add/subtract operations per vector item for assessing the error of rounding, so do 4 operations per data item. Thus, theoretical peak for twofold-fast method would be 4 times lower, 6 gigaflops for float and 3 gigaflops for double type. Peak for Kahan summation is the same, as it also does 3 extra operations per data item. Look at results for psumtf, pdottf, psumk, pdotk.

More rigorous twofold summation takes 7 add/subtract operations per data item, so peak performance for psumt, pdott functions is 3.4 (≈24/7) gigaflops for floats, and 1.7 (≈12/7) gigaflops for doubles. Look at the table for the results of the psumt, pdott functions.

Test results for psumd, pdotd which use double-precision accumulator for float data is 12 gigaflops, the theoretical peak for summation with double precision.

Here we do not test quad precision, which AVX doesn't support.

To conclude: the p-functions look optimized very high and utilize processor capacity at theoretical peak.

## Highly optimized functions

Highly optimized function names start with "h" prefix like hsum, hdot, etc. These functions share same code with p-prefixed functions psum, pdot, etc. The difference is that h-functions do read memory.

Please look at the testing results with small data arrays that can fit into processor fastest L1 cache:

Reading does not affect performance of the hsum function, so it shows the peak 24 gigaflops for float and 12 gigaflops for double. For hdot however, reading impact is substantial, as reading is apparently limited like demonstrated with the read2 test. So hdot performs twice slower, at 12 gigaflops for single precision and 6 gigaflops for double.

Reading from L1 cache does not limit twofold and Kahan methods as they more depend on processor capacity. Reading from L1 appears to impact hsumd, hdotd; unexpectedly for me.

If we look at results with large arrays that cannot fit into processor cache, performance of simple hsum, hdot looks completely determined by the reading bandwidth in this case. Compare hsum, hdot results with read1, read2 for large arrays.

Dependency for twofold and Kahan methods look more complicated. Please note however, that twofold and Kahan methods perform at nearly same level as simple summation for large arrays, because of the negative impact of memory reading. Processor is underused while fetching data, and twofold and Kahan methods can effectively leverage this resource for assessing or compensating rounding errors.

Note, that twofold fast summation for double precision data works 30-100 times faster than if we used quad precision for improving accuracy. Look at sumd, dotd results (quad is g++ __float128 here).

## Naïve vectorization for AVX

For making a wider picture, the table includes v-prefixed functions vsum, vdot, vsumd, vdotd compiled for AVX with fast-math and maximal optimization options:

      g++ –ffast-math –O3 –mavx …

The other a-prefixed functions asum, adot, etc. are my hand-made made naïve vectorization for AVX. The code is pretty much the same as for v-functions with two important exceptions: a-functions do not accept unaligned data, and I vectorized trickier twofold and Kahan methods that compiler cannot.

Both hand-made simple asum, adot and vsum, vdot perform at around at 1/3 of processor peak, 8 gigaflops for float and 4 gigaflops for double type. The point here is unresolved data dependency. Next AVX instruction cannot start before previous ends, as next instruction uses result of the previous. For the tested processor, latency of AVX summation is 3 ticks, so we observe 3 times slower performance.

The h-prefixed functions we consider above feed the conveyor every tick, so gain 100% of CPU capacity.

## Scalar u-functions

Next group is functions with names prefixed with "u", usum, udot, etc. These functions coded purely in C++, and compiled with best optimizations but with strict math option and without vectoring for AVX. The coding is tricky enough to resolve the data dependency, so enable nearly peak performance for non-vectored case, which is one summation per processor tick, or 3 gigaflops if doing at 3 gigahertz.

Twofold fast summation performs at 750 megaflops with such coding, as expected because it does 4 of add/subtract operations per data item. Twofold rigorous method performs at nearly 430 megaflops, which is expected 1/7 of processor peak, as it does 7 operations per data item.

Kahan compensatory summation with 4 operations per item performs at 750 megaflops as expected.

## Scalar functions

The last but not least group is functions with no prefix in the name, sum, dot, etc. This is straightforward basic variant of coding, compiled with maximal optimization but with strict math. Note, that we cannot compile twofold and Kahan methods with fast-math option, as associative optimization can damage the tricks with rounding errors on which these methods rely.

Performance of all methods looks disappointing in this test, but even so twofold-fast summation for double precision is still 20 times faster than using quad precision for improving or assessing accuracy.

## Conclusion

SIMD vectorization is the must for enjoying twofold bonuses with minimal performance loss. With AVX vectorization, modern processor is faster than memory, so twofold can utilize underused CPU capacity for additional useful job, improving accuracy and automating control of errors accumulation.

This seems a general trend. Due to more cores and larger SIMD conveyors, math software can typically use only part, maybe 10% of processor capacity. Leveraging the remaining 90% gap is the opportunity; and twofold fast or Kahan summation is an obvious way for utilizing this gap.

## Algorithms

Twofold approach combines of well-known techniques investigated in 1960 and 1970. Specifically, twofold bases on two theorems usually credited to Dekker and Knuth. Following is citation from the Shewchuk paper [4] dedicated to automatically adapting accuracy to specific computations. Here, symbols ⊕ and ⊖ mean floating-point add/subtract operations with correctly rounding to even:

> Theorem 1 (Dekker). Let $a$ and $b$ be $p$-bit floating-point numbers such that $|a| \geq |b|$. Then the following algorithm will produce a non-overlapping expansion $x + y$ such that $a + b = x + y$, where $x$ is approximation to $a + b$ and $y$ represents the round-off error in the calculation of $x$.
>
> FAST-TWO-SUM($a$, $b$)
> 1. $x \leftarrow a \oplus b$
> 2. $b' \leftarrow x \ominus a$
> 3. $y \leftarrow b \ominus b'$
> 4. return ($x$, $y$)
>
> Theorem 2 (Knuth). Let $a$ and $b$ be $p$-bit floating-point numbers, where $p \geq 3$. Then the following algorithm will produce a non-overlapping expansion $x + y$ such that $a + b = x + y$.
>
> TWO-SUM($a$, $b$)
> 1. $x \leftarrow a \oplus b$
> 2. $b' \leftarrow x \ominus a$
> 3. $a' \leftarrow x \ominus b'$
> 4. $b^\# \leftarrow b \ominus b'$
> 5. $a^\# \leftarrow a \ominus a'$
> 6. $y \leftarrow a^\# \oplus b^\#$
> 7. return ($x$, $y$)

The 2$^{nd}$ theorem eliminates the need in checking if $|a| \geq |b|$ for the cost of 3 additional operations. With modern processors, 3 additional add/subtract operations would cost less than conditional branching, so the software like for double-double and quad-double arithmetic use to base on the Knuth theorem.

Twofold fast summation however bases on the Dekker theorem, as for the goals of twofold summation, there is no much practical benefit in the more rigorous formula. Following is non-optimized code:

> Example 1: Twofold type definition and fast summation algorithm
>
> ```
> // Define twofold<number>, assume number is float or double
> template<typename number> struct twofold { number value, error; };
>
> // Twofold fast summation (based on Dekker theorem)
> template<typename number> twofold<number> sumtf(int m, number data[]) {
>     number s=0, e=0;
>     for (int i=0; i<m; i++) {
>         number y, t, c;
>         y = data[i];
>         t =  s + y;
>         c = (t - s) - y;
>         e =  e - c;
>         s =  t;
>     }
>     twofold<number> result;
>     result.value = s;
>     result.error = e;
>     return result;
> }
> ```

Such code is very similar to the Kahan compensatory summation, as it follows from the Wikipedia:

Example 2: Kahan compensatory summation

```cpp
// Kahan compensatory summation, number is float or double
//   http://en.wikipedia.org/wiki/Kahan_summation_algorithm
template<typename number> number sumk(int m, number data[]) {
    number s=0, c=0;
    for (int i=0; i<m; i++) {
        number y, t;
        y = data[i] - c;
        t =  s + y;
        c = (t - s) - y;
        s =  t;
    }
    return s;
}
```

The difference is the accent: while Kahan method strives minimizing the rounding error, twofold fast summation computes the *value* exactly like a simple direct summation would and tries estimating its inaccuracy by collecting the round-offs with the separate *error* variable. Provided the round-offs are much smaller than the data items, the *error* must reasonably estimate inaccuracy of the *value*.

Once we get an *error* estimate, what shall we do with it? We can just ignore it, or use it for improving accuracy. The sum of *value* + *error* must be more accurate than *value* alone, even if evaluated with the number's original precision. If we evaluate *value* + *error* with higher precision, my result must be more accurate than for Kahan summation, because the *error* can accumulate more of the significant bits.

But I think the better idea would be to interrupt computations if *error* gets too large. Such situation would basically signal that the used floating-point precision looks not enough for the specific task. Maybe it is worth to think of using higher-precision numbers, like use double instead of float.

Here is the code (non-optimized variant) for the Kahan and simple summation functions:

Example 3. Twofold rigorous and direct summation

```cpp
// Twofold summation (rigorous, based on Knuth theorem)
template<typename number> twofold<number> sumt(int m, number data[]) {
    number s = 0, e = 0;
    for (int i = 0; i<m; i++) {
        number t, y, yt, dy, ds;
        t  = s + (y  = data[i]);
        dy = y - (yt = t - s);
        ds = s - (t - yt);
        e += ds + dy;
        s  = t;
    }
    twofold<number> result;
    result.value = s;
    result.error = e;
    return result;
}

// Simple direct summation, assume number is float or double
template<typename number> number sum(int m, number data[]) {
    number s=0;
    for (int i=0; i<m; i++)
        s += data[i];
    return s;
}
```

```
// Summation with double accumulator
double sumd(int m, float data[]) {
    double s=0;
    for (int i=0; i<m; i++)
        s += data[i];
    return s;
}
```

Similar dot product functions would do multiplication like `a[i]*b[i]` in place of `data[i]`.

Sources for the other functions used in this article are available at my Web page [7].

## Accuracy

Following table shows results of accuracy testing for different summation methods:

Table 2. Random numbers summation accuracy (relative error to precise solution)

| Type | Random generator | Interval | Direct | Kahan | Twofold fast |
|---|---|---|---|---|---|
| float | Numerical Recipes | [0,1] | $-1.07759_{10}$-07 | $1.72372_{10}$-08 | $-5.24003_{10}$-11 |
| float | Numerical Recipes | [-1,1] | -0.000273821 | $3.69634_{10}$-09 | 0 |
| float | MMIX by D. Knuth | [0,1] | $-4.74619_{10}$-06 | $5.19367_{10}$-09 | $9.09826_{10}$-11 |
| float | MMIX by D. Knuth | [-1,1] | $1.65869_{10}$-06 | $-2.07055_{10}$-08 | 0 |
| double | Numerical Recipes | [0,1] | $2.72008_{10}$-17 | $2.72008_{10}$-17 | 0 |
| double | Numerical Recipes | [-1,1] | $3.20017_{10}$-13 | $-1.06717_{10}$-16 | 0 |
| double | MMIX by D. Knuth | [0,1] | $-1.63303_{10}$-14 | $-2.73781_{10}$-17 | 0 |
| double | MMIX by D. Knuth | [-1,1] | $5.17981_{10}$-14 | $-1.14555_{10}$-17 | 0 |

The test creates data array of 1-million elements fulfilled with random numbers uniformly distributed in the interval $[0,1]$ or $[-1,1]$. The test sums the data with three methods, directly, with Kahan method, and with twofold fast method. Simultaneously, the test sums the same data with quad precision, and considers this extra-precise sum as a good approximation to the exact solution. The table displays the relative error of the tested methods against the quad-precision solution. The test uses random number generators from Wikipedia: http://en.wikipedia.org/wiki/Linear_congruential_generator

Because most of the test array elements are of same order by magnitude, this must be a good case for twofold summation. Indeed, twofold fast sum appears same accurate as quad-precision sum in most of the test cases here.

Now let us see how twofold techniques works in the "100 hours" test which I know from Marius Cornea, and which I think must be a well-known test:

Imagine a timer that ticks every 1/10 of second and on each tick adds 0.1 to a single-precision counter. What time would such clock show after 100 hours? With single precision, round-off errors accumulate too much, so summation shows ~96 hours instead of 100. Let us see how various summation methods would address this problem. Here "double" method is same timer but with double precision counter:

Table 2. Results and deviation in "100 hours" test

| Method | Result | Deviation | Estimate |
|---|---|---|---|
| Direct | 96.3958 | 3.60423 | – |
| Double | 100 | $1.49012_{10}$-06 | – |
| Kahan | 100 | 0 | – |
| Twofold | 99.9359 | 0.0641498 | 3.54008 |

For twofold method, this table shows *value + error* sum in the Result column, and shows the *error* in the Estimate column. In this test, the *value* deviation from exact solution is much higher than the estimate $O(\sqrt{N}\varepsilon) \approx 1.13_{10}$-4 (0.0113 hours). The *value + error* sum also deviates much wider than the estimate $O(N\varepsilon^2) \approx 1.28_{10}$-6 hours, but appears enough to signal that single precision is not appropriate here.

Now let us deduce the asymptotic estimate for accuracy of twofold summation:

For direct summation of $N$ numbers, best-case inaccuracy estimate is zero, in such improbable case if all calculations were exact. Worst-case estimate for relative error is $O(N\varepsilon)$, if we lose 1/2 bit per each $x_n$; here $\varepsilon$ is $\frac{1}{2}\text{ulp}(1)$ for the used floating-point precision (ulp is "unit in last position"). As the average-case, it is useful to assume, that the round-offs might have different signs. Assuming sign is random, average-case estimate for relative error in direct summation would be $O(\sqrt{N}\varepsilon)$.

Now consider average-case estimate for inaccuracy of *value* + *error* sum of twofold result:

Relative inaccuracy of the *value* would be $O(\sqrt{N}\varepsilon)$, or in other words, the *value* deviates from the exact result approximately by *value*·$O(\sqrt{N}\varepsilon)$. Similarly, absolute deviation between the *error* and the exact of the accumulated error is *error*·$O(\sqrt{N}\varepsilon)$. Presumably, *error* estimates the deviation of the *value*, so let us suppose that the *error* approximately equals it, *error* ≈ *value*·$O(\sqrt{N}\varepsilon)$. By substituting this into previous formula, we can estimate deviation of *error* like (*value*·$O(\sqrt{N}\varepsilon)$)·$O(\sqrt{N}\varepsilon)$, or *value*·$O(N\varepsilon^2)$.

Provided *value* + *error* approximately equals to *value*, its relative inaccuracy is $O(N\varepsilon^2)$ in average case.

Such estimate must be better than $O(\sqrt{N}\varepsilon)$ for direct sum, if $N$ is not too large. Actually, just better is not enough; we want *value* + *error* be much more accurate than *value* alone, so *error* can estimate the deviation of *value*. That is, we need $O(N\varepsilon^2)$ be much less than $O(\sqrt{N}\varepsilon)$. It is possible if $N$ is much less than $1/\varepsilon^2$. I guess, twofold must work fine if $N < 1/\varepsilon$, that is $N$ must not exceed $2^{24} \approx 10^7$ for single precision and $2^{53} \approx 10^{16}$ for double.

# Conclusion

Floating-point arithmetic works perfectly in millions of applications, but fails sometimes due to rounding errors accumulation. Debugging such failures is difficult; and ideally, computer should resolve such cases automatically or at least track it and signal clearly if a failure. Obvious idea for that is simply increasing precision of floating-point numbers, which we often do by using double precision instead of single.

Obvious next step would be using quad precision, if double appears not enough. This is good approach, except available implementations of quad arithmetic are very slow, 100 or even more times slower than double precision with modern processors. Double-double and many-of-doubles techniques allow better precision in a different manner.

Tracking round-offs was one of the reasons for introducing interval arithmetic in $1950^{th}$ or even earlier. However, intervals look more appropriate for their main reason, guaranteeing the boundaries for result. Such more ambitious goal makes intervals to get much wider than we need for tracking rounding errors in an "average case" assumption, as interval summation must cover all cases including the worst one.

Twofold fast summation combines the double-double and interval techniques and represent an exact real with *value + error* pair of floating-point numbers. Similar to double-double, *value + error* must be more accurate than *value* alone, so *error* can estimate interval between the *value* and its exact target. Unlike double-double, twofold strives to estimate inaccuracy rather than compensate it.

Unlike intervals, twofold estimate cannot guarantee the boundaries for the exact result, so can provide completely wrong perception of the result and of the accuracy of the approximation. But I think, more probably twofold would more-or-less correctly estimate accumulation of inaccuracy in computations. Improving the actual precision with *value + error* would be a free bonus.

Key advantage of twofold approach is its high performance with modern processors. If leveraging SIMD vectorization appropriately, twofold would work only 1-4 times slower than the best of top-performing manually optimized direct summation. Anyway, even in worst case twofold assessment works at least 20-100 times faster than if we used quad precision for tracking inaccuracy.

Moreover, in a typical situation if at least part of data is not pre-loaded to processor cache, twofold fast summation can perform at the same level as highly optimized direct method. Because processor is much faster than memory, direct method can use only part of CPU capacity, maybe 10%. Twofold fast method can utilize the remaining 90% and do the useful extra job while fetching the data from memory.

Underusing CPU capacity would be the general trend, as newer processors with more cores and wider SIMD conveyor would increase the performance gap with memory. It becomes profitable utilizing this gap for something useful; and twofold fast summation is an obvious example.